\documentclass[12pt]{article}
\usepackage{amssymb, amsmath, amscd}
\begin{document}
\renewcommand{\theequation}{\arabic{section}.\arabic{equation}}
\newtheorem{theorem}{Theorem}[section]
\newtheorem{lemma}{Lemma}[section]
\newtheorem{definition}{Definition}
\newtheorem{pro}{Proposition}[section]
\newtheorem{cor}{Corollary}[section]
\newtheorem{remark}{Remark}[section]
\newcommand{\n}{\nonumber}
\newcommand{\tv}{\tilde{v}}
\newcommand{\tw}{\tilde{\omega}}
\renewcommand{\t}{\theta}
\newcommand{\w}{\omega}
\newcommand{\abs}[1]{\left| #1 \right|}
\newcommand{\bke}[1]{\left( #1 \right)}
\newcommand{\calC}{{ \mathcal C  }}
\newcommand{\calS}{{ \mathcal S  }}
\newcommand{\calP}{{ \mathcal P  }}
\newcommand{\calJ}{{ \mathcal J  }}
\newcommand{\calI}{{ \mathcal I }}
\newcommand{\R}{\Bbb R}
\newcommand{\e}{\varepsilon}
\renewcommand{\a}{\alpha}
\renewcommand{\l}{\lambda}
\newcommand{\vare}{\varepsilon}
\newcommand{\s}{\sigma}
\renewcommand{\o}{\omega}
\renewcommand{\O}{\Omega}
\newcommand{\bb}{\begin{equation}}
\newcommand{\ee}{\end{equation}}
\newcommand{\bq}{\begin{eqnarray}}
\newcommand{\eq}{\end{eqnarray}}
\newcommand{\bqn}{\begin{eqnarray*}}
\newcommand{\eqn}{\end{eqnarray*}}
\title{Remarks on regularity conditions of the  Navier-Stokes equations}
\author{Dongho Chae\thanks{This work was supported partially by  KRF Grant(MOEHRD, Basic
Research Promotion Fund) and the KOSEF Grant no.
R01-2005-000-10077-0.}\\
Department of Mathematics\\
              Sungkyunkwan University\\
               Suwon 440-746, Korea\\
              {\it e-mail : chae@skku.edu}}
 \date{}
\maketitle
\begin{abstract}
Let $v$ and $\o$ be the velocity and the vorticity of the a suitable
weak solution of the 3D Navier-Stokes equations in a space-time
domain containing $z_0 =(x_0, t_0)$, and $Q_{z_0, r} =B_{x_0,
r}\times (t_0-r^2 , t_0)$  be a parabolic cylinder in the domain. We
show that if $v\times \frac{\o}{|\o|}\in L^{\gamma, \alpha}_{x,t}
(Q_{z_0, r})$ or $\o \times \frac{v}{|v|}\in L^{\gamma,
\alpha}_{x,t} (Q_{z_0, r})$, where $L^{\gamma, \alpha}_{x,t}$
denotes the Serrin type of class, then $z_0$ is a regular point for
$v$. This refines previous  local regularity criteria for the
suitable weak solutions.
\end{abstract}
\section{Introduction}
 \setcounter{equation}{0}

The Navier-Stokes equations in a domain $\Omega \in \Bbb R^3$ are
the following.
$$
(NS)\left\{ \aligned &\frac{\partial v}{\partial t} + (v \cdot
\nabla ) v=-\nabla p +
 \Delta v, \qquad (x,t)\in \Omega\times (0, T)\\
&\mathrm{div}\, v=0,\qquad (x,t)\in \Omega\times (0, T)\\
 & v(x,0)=v_0(x), \qquad  x\in \Omega\\
\endaligned\right.
$$
where $v=(v_1, v_2, v_3 )$, $v_j =v_j (x, t)$, $j=1,2,3$, is the
velocity of the flow, $p=p(x,t)$ is the scalar pressure, and $v_0$
is the given initial velocity satisfying div $v_0=0$. The global in
time existence of a smooth solution  to the system (NS) is an
outstanding open problem in mathematics, and is chosen as one of the
seven millennium problems by Clay Institute. One traditional
approach to the problem is to prove global in time existence of weak
solutions, and the prove their regularity.  A notion of weak
solution of (NS) was introduced, and its global in time existence in
$\Bbb R^3$ was proved by Leray in \cite{ler}. Later, Hopf proved
existence of weak solution in a bounded domain in \cite{hop}. After
that there are numerous conditional regularity results on the weak
solutions, imposing integrability conditions on the velocity or the
vorticity, which guarantees regularity of the weak solutions(see
e.g. \cite{ser, ohy, pro, lad, fab, str, bei1,bei2,bei3,tia, esc,
neu, cha1,cha2,cha3}). For the local analysis of the regularity
properties of weak solutions Caffarelli-Kohn-Nirenberg introduced
the notion of suitable weak solutions and proved its partial
regularity  as well as global in time existence(\cite{caf}). A
refined definition of suitable weak solutions, using a stronger
condition for pressure, which we adopt here, was introduced by Lin
in \cite{lin}. Let $Q_T=\Omega\times (0,T)$.
 For a point $z=(x,t)\in Q_T$, we denote below
  \bqn B_{x,r}= \{ y\in \Bbb R^3:
|y-x|<r \}, \quad Q_{z,r} =B_{x,r} \times (t-r^2, t).
 \eqn
\begin{definition}
 A pair $(v, p)$ of measurable functions is a suitable weak solution
 of (NS) if the
following conditions are satisfied:
\begin{itemize}
\item[(i)]
$
v\in L^{\infty}(0,T;L^2(\Omega))\cap L^2(0,T;W^{1,2}(\Omega)),\quad
p \in L^{\frac32}(Q_T).
$
\item[(ii)] The following integral identity holds
$$\int_{Q_T} \left[ -v\cdot \partial_t \varphi +(v\cdot \nabla )v \cdot
\varphi +\nabla v : \nabla \varphi \right]dxdt=\int_{\Omega} v_0
\cdot \varphi (x, 0)dx
$$
for all vector test functions $\varphi \in [C_0 ^\infty (\Omega
\times [0, T))]^3$.
\item[(iii)] The pair $(v,p)$ satisfies the local energy
inequality,
\[
\int_{\Omega}\abs{v(x,t)}^2\phi(x,t) dx
+2\int_0^t\int_{\Omega}\abs{\nabla v(x,\tau)}^2\phi(x,\tau) dxd\tau
\]
$$
\leq\int_0^t\int_{\Omega}\bke{\abs{v}^2(\partial_t\phi+\Delta\phi)
+(\abs{v}^2+2p)v\cdot\nabla\phi }dxd\tau
$$
 for almost all $t\in
(0,T)$ and for all nonnegative scalar test function $\phi\in
C_0^{\infty} (Q_T)$.
\end{itemize}
\end{definition}
We say that a weak solution $v$ is  regular at $z$, if $v$ is
bounded in $Q_{z,r}$ for some $r>0$. Such point $z$ is called a
regular point. A point in $Q_T$, which is not regular, is called a
singular point. Caffarelli-Kohn-Nirenberg showed  that the one
dimensional  Hausdorff measure of the set $\calS$ of possible
interior singular points of suitable weak solutions is
zero(\cite{caf}), which refines the previous results due to
Scheffer(\cite{sch}).

In this paper our aim is to obtain refined versions of regularity
conditions for velocity and vorticity for suitable weak solutions,
incorporating the directions of each vector field as well as the
magnitudes. Our conditions are not directly on the velocity or
vorticity, but on the orthogonal component of velocity to vorticity
direction, or on the orthogonal component of vorticity to velocity
direction. The associated integral norms are scaling invariant. More
precisely our main theorem is the following.

\begin{theorem}
Let $z_0=(x_0,t_0) \in Q_{T}$ with  $\bar{Q}_{z_0,r}\subset Q_{T}$,
and  $(v, p)$ be a suitable  weak solution of (NS) in $Q_T$ with the
vorticity $\o=$curl$\, v$, where the derivatives are in the sense of
distribution. Suppose $v$ and $\o$ satisfy one of the following
conditions:
\begin{itemize}
\item[(i)] The following inequality holds true,
 \bb\label{th0}
 \left\|v\times \frac{\o}{|\o|}\right\|_{L^{3, \infty}_{x,t}
 (Q_{z_0,r})}\leq \vare_0
 \ee
 for an absolute constant $\vare_0 >0$ to be determined in the proof
 below.
\item[(ii)] There exist $\gamma \in (3, \infty]$, $\alpha \in [2,
\infty)$ with $3/\gamma +2/\alpha \leq 1$ such that
  \bb\label{th1}
  v\times \frac{\o}{|\o|} \in L^{\gamma, \alpha}_{x,t} (Q_{z_0,r}).
  \ee
 \item[(iii)] There exist $\gamma \in [2, 3]$, $\alpha \in [2,4]$
 with $3/\gamma +2/\alpha \leq 2$ such that
 \bb\label{th2}
 \o \times \frac{v}{|v|}\in L^{\gamma, \alpha}_{x,t} (Q_{z_0,r}).
  \ee
  Then,  $z_0$ is a regular point.
  \end{itemize}
\end{theorem}
\begin{remark}
In a recent preprint J. Wolf(\cite{wol}) proved an
$\vare-$regularity type of theorem related to (iii) of the above
theorem as follows: There exists $\vare_0 ^*$ such that if a
suitable weak solution in $Q_{z_0,r}$ satisfies
$$
\lim\sup_{\rho\to 0+} \frac{1}{\rho} \int_{Q_{z_0, \rho }} \left|\o
\times \frac{v}{|v|} \right|^2 dxdt \leq \vare_0 ^* ,
$$
then, $z_0$ is a regular point.  If the condition (iii) holds true,
then by H\"{o}lder's inequality and the absolute continuity of the
integral norms we have
$$
\lim\sup_{\rho\to 0+}\frac{1}{\rho} \int_{Q_{z_0, \rho }} \left|\o
\times \frac{v}{|v|}\right|^2 dxdt \leq \lim\sup_{\rho\to 0+}\rho
^{2(2-\frac{3}{\gamma}-\frac{2}{\alpha})} \left\|\o \times
\frac{v}{|v|} \right\|_{L^{\gamma, \alpha}_{x,t} (Q_{z_0, \rho})}=
0.
$$
Hence, the part (iii) of the above theorem is implied by the main
theorem of \cite{wol}. The part (i), (ii) of the above theorem,
however, have no direct implication relationships with that result.
Moreover, the proof of part (iii) of the above theorem  given below
is much simpler than that of \cite{wol}.
\end{remark}

\begin{remark}
 We say $v$ is a
Beltrami flow in $Q_{z_0,r}$ if  $v\times \o=0$ in $Q_{z_0,r}$. In
the study of physics of turbulent flows the Beltrami structure has
important roles(see e.g. \cite{con2, pel} and the references
therein).
 The condition that $v\times \frac{\o}{|\o|}$ or $\o \times
\frac{v}{|v|}$  is controllable in a space-time region  implies
intuitively that the weak solutions are not far from the Beltrami
flows in that region in an appropriate sense, and the above theorem
says that this implies regularity of the flows in that region.
\end{remark}

\section{Proof of Theorem 1.1}
 \setcounter{equation}{0}
 Before starting our proof we recall previous results concerning the
 notion of {\it{an epoch of
possible irregularity}} of the weak solution of the Navier-Stokes
equations.  It is known that for weak solutions there exists a set
$E\subset I=[0,T]$ such that $E$ is closed, of $1/2$-dimensional
Hausdorff measure zero, and solutions are regular in $I\setminus E$
(\cite{ler, gal, foi}).  Moreover, the set $E$ can be written as
$I\setminus\cup_{i\in\calJ}I_i$, where set $\calJ$ is at most
countable, and $I_i=(\alpha_i,\beta_i)$ are disjoint open intervals
in $[0,T]$. Following \cite{gal}, we call the instant time $\beta_i$
 an epoch of possible irregularity.  We recall a fact proved by
Neustupa and Penel in \cite{neu} on  the epoch of possible
irregularity  for suitable weak solutions.
\begin{lemma}\label{lem1}
Let $z_0=(x_0,t_0)\in Q_T$. Suppose $v$ is a suitable weak solution
of (NS)  in $Q_{T}$ and $t_0$ be an epoch of possible irregularity.
Then there exist positive numbers $\tau$, $r_1$, and $r_2$ with
$r_1<r_2$ such that the followings are satisfied:
\begin{itemize}
\item[(a)]
$\tau$ is sufficiently small so that $t_0$ is only one epoch of
possible irregularity in time interval $[t_0-\tau,t_0]$.
\item[(b)]
The closure of $B_{x_0, r_2} \times (t_0-\tau, t_0)$ is contained in
$Q_T$, i.e. $\overline{B_{x_0, r_2}} \times [t_0-\tau, t_0] \subset
Q_T$.
\item[(c)] $\bke{ (\overline{B_{x_0, r_2}}-B_{x_0, r_1})\times [t_0-\tau,
t_0]}\cap \calS= \phi$, where $\calS$ is the set of possible
singular points  of $v$.
\item[(d)]
 $v$, $v_t$, and $p$ are, together with
all their space derivatives, continuous on $(\overline{B_{x_0,
r_2}}-B_{x_0, r_1})\times [t_0-\tau, t_0]$.
\end{itemize}
\end{lemma}
Next we recall the following result proved in \cite{caf}, a
corollary of which will be used in the proof of our main theorem.
\begin{pro}
There exists an absolute constant $\vare_1>0$ with the following
property. If $(v,p)$ is a suitable weak solution of (NS) near $z_0$
and if
 \bb\label{pro1} \lim\sup_{\rho \to 0+} \frac{1}{\rho}\int_{Q_{z_0, \rho}} |\nabla v
 |^2 dxdt \leq \vare_1,
 \ee
 then $z_0$ is a regular point.
\end{pro}
As an immediate corollary we have the following local regularity
criterion, which is a local version of the one obtained in
\cite{bei2}.
\begin{cor}
If $(v,p)$ is a suitable weak solution of (NS) near $z_0$, and if
 either
 $$ \|\nabla v \|_{L^{\frac32, \infty} _{x,t} (Q_{z_0, r})} \leq
 \vare_1 ,
 $$
 where $\vare_1$ is the constant in Proposition 2.1, or there exist
 $\gamma \in (3/2, \infty]$ and $\alpha\in [2, \infty)$
with $3/\gamma +2/\alpha \leq 2$ such that
$$ \nabla v \in L^{\gamma, \alpha} _{x,t} (Q_{z_0, r}), $$
then, $z_0$ is a regular point.
\end{cor}
\noindent{\bf Proof.} Similarly to Remark 1.1 we observe
 \bqn \lim\sup_{\rho\to 0+}\frac{1}{\rho}
\int_{Q_{z_0, \rho }} \left|\nabla v\right|^2 dxdt &\leq
\lim\sup_{\rho\to 0+}\rho ^{2(2-\frac{3}{\gamma}-\frac{2}{\alpha})}
\left\|\nabla v\right\|_{L^{\gamma, \alpha}_{x,t} (Q_{z_0, \rho})}\\
&\left\{
\aligned = 0 \, \quad&\mbox{if} \quad \, \mbox{$\gamma >3/2$ and
$3/\gamma +2/\alpha \leq 2$}\\
\leq \vare_1  \, \quad&\mbox{if} \,\quad (\gamma,
\alpha)=(3/2,\infty)
\endaligned \right.
\eqn
 by the
H\"{o}lder inequality. Then the conclusion is immediate by
Proposition 2.1. $\square$\\
\ \\
 \noindent{\bf Proof of Theorem 1.1}
We first assume that $t_0$  is an epoch of possible irregularity for
$v$ in $Q_{z_0, r}$. Suppose  that $0<r_1 <r_2 <r$, and $r ^2 <
\tau$ are the positive numbers in Lemma \ref{lem1}. Below, we denote
$B_1=B_{x_0, r_1}$ and $B_2=B_{x_0,r_2}$.  Following \cite{neu}, we
choose a cut-off function $\varphi\in C^{\infty}_0(B_2)$ such that
$\varphi = 1$ on $B_1$, and set $u=\varphi v-V$, where $V\in C^2_0
(B_2\backslash \overline{B_1})$ satisfies $\mbox{div } V= v\cdot
\nabla \varphi$.
 In particular, all the spatial derivatives of $V$ and $\frac{\partial V}{\partial t}
$ are smooth. Using the well-known form of the Navier-Stokes
equations,
$$
\frac{\partial v}{\partial t} -v\times\o=-\nabla  (p+\frac12 |v|^2)+
\Delta v,
$$
 one can check easily
that $u$  satisfies the following equations:
 \bb\label{eqforu}
 \frac{\partial
u}{\partial t}-\varphi v\times\o=h-\nabla \left( \varphi (p+\frac12
|v|^2)\right)+ \Delta u, \quad\mbox{div }u=0,
 \ee
  where we set
 \bqn
  h&=&-\frac{\partial V}{\partial t}+(p+\frac12 |v|^2 )\nabla \varphi- v \Delta \varphi
  -2 (\nabla \varphi \cdot\nabla )v
 + \Delta V .
 \eqn
   We observe that $h(\cdot,t)$ is supported on
$\bke{ \overline{B_{2}}\setminus B_{1}}$ for each $t\in [t_0-\tau,
t_0)$, which is sufficiently smooth in the region.  Operating $D$ on
(\ref{eqforu}), and taking $L^2 (B_{x_0, r_2}) $ inner  product it
by $Du$, we obtain, after integration by part
  \bq\label{eqfordu}
  \lefteqn{\frac12 \frac{d}{dt} \|Du \|_{L^2} ^2 + \|D^2 u\|_{L^2} ^2 =
  -(\varphi v\times \o , D^2 u )_{L^2} -(D^2 u, h) }\hspace{.1in}\n \\
  && \leq |(\varphi v\times \o , D^2 u )_{L^2} | + \frac{1}{8} \|D^2 u\|_{L^2} ^2
  +C \|h\|_{L^2}^2 ,
  \eq
  where(and below) we used simplified notation for the $L^p-$norm in $B_2$,
  $$\|f\|_{L^p} =\|f\|_{L^p (B_2)}, \quad p\in [1, \infty],$$
  unless other domain is specified.
  Let us set $\xi =\o/|{\o} |$.
 We  estimate the nonlinear term  as follows:
 \bq\label{estfornon}
|(\varphi v \times \o , D^2 u)_{L^2}|&\leq&
\int_{B_2} |v\times \xi|| \varphi \o | |D^2 u |dx\n \\
&\leq& \int_{B_2} |v\times \xi ||\varphi Dv||D^2 u|dx\n \\
 &=&\int_{B_2} |v\times \xi | |Du -v\nabla \varphi +DV||D^2 u |dx\n \\
&\leq& \int_{B_2} |v\times \xi | |Du||D^2 u |dx +\int_{B_2} |v\times
\xi ||g| |D^2 u|dx \n \\
&=&I_1 +I_2,
  \eq
  where we set $g=v\nabla \varphi -DV.$
  Since $g$ is a smooth function supported on  $(B_2 \setminus
\bar{B_1}) \times (t_0- \tau ,t_0 ]$, we estimate $I_2$ simply as
 \bq\label{estfori2}
 I_2 \leq \|g\|_{L^\infty}\|v\|_{L^2} \|D^2 u\|_{L^2}\leq
 C\|v\|_{L^2}^2 +\frac{1}{4} \|D^2 u\|_{L^2}^2.
 \eq
  We first assume the condition (i) of Theorem 1.1 holds true. In this case we estimate
  \bb\label{l3}
   I_1\leq \|v\times \xi \|_{L^3} \|Du\|_{L^6} \|D^2 u\|_{L^2}\leq C\|v\times \xi \|_{L^3}
  \|D^2 u\|_{L^2}^2.
  \ee
Combining estimates (\ref{eqfordu})-(\ref{l3}) together, we have
 \bq\label{integ}
\frac{d}{dt} \|Du \|_{L^2} ^2 + \|D^2 u\|_{L^2} ^2 &\leq&
C_1\|v\times \xi \|_{L^3}\|D^2 u\|_{L^2}^2 +C
\|h\|_{L^2}^2+ C\|v\|_{L^2}^2\n \\
&\leq & C_1\vare_0\|D^2 u\|_{L^2}^2 +C \|h\|_{L^2}^2+ C\|v\|_{L^2}^2
 \eq
for $t\in (t_0 -r_2 ^2, t_0]$, and for an absolute constant $C_1$.
 If $C_1 \vare_0 < 1$, then integrating  (\ref{integ}) in time over $[t_0-r_2^2,t_0]$, we can obtain
\bqn
  \lefteqn{\sup_{t_0 -r_2^2<t<t_0}\|Du (\cdot ,t)\|_{L^2 } ^2 \leq  \|Du (\cdot
,t_0-r_2 ^2)\|_{L^2 }
 ^2 +C\int_{t_0 -r_2 ^2} ^{t_0} \|v\|_{L^2} ^2 dt}\hspace{1.5in}\n \\
 && +C\int_{t_0 -r_2 ^2} ^{t_0} \|h\|_{L^2} ^2 dt
 <\infty.
 \eqn
Hence, $Du\in L^{2, \infty}_{x,t} (Q_{z_0 ,r_2}) $, and  therefore
$Dv\in L^{2, \infty}_{x,t} (Q_{z_0 ,r_1})$. Applying Corollary 2.1,
we conclude that $z_0$ is a regular point.  Next, we assume that the
condition (ii) of Theorem 1.1 holds true, and estimate
 \bq\label{estfori1}
 I_1&\leq& \|v\times \xi \|_{L^\gamma} \|Du
 \|_{L^{\frac{2\gamma}{\gamma -2}}} \|D^2 u \|_{L^2} \n \\
 &\leq & C\|v\times \xi \|_{L^\gamma}\|Du
 \|^{1-\frac{3}{\gamma}}\|D^2 u  \|_{L^2} ^{1+\frac{3}{\gamma}}\n \\
 &\leq & C\|v\times \xi \|_{L^\gamma}\|Du
 \|^{1-\frac{3}{\gamma}}\|D^2 u \|_{L^2} ^{1+\frac{3}{\gamma}}\n \\
 &\leq& C\|v\times \xi \|_{L^\gamma} ^{\frac{2\gamma}{\gamma -3}}
 \|Du
 \|_{L^2} ^2 +\frac{1}{4} \|D^2 u \|_{L^2} ^2,
 \eq
 where we used the interpolation inequality,
  $$ \|Du
 \|_{L^{\frac{2\gamma}{\gamma -2}}} \leq C \|Du
 \|^{1-\frac{3}{\gamma}} _{L^2} \|D^2 u\|_{L^2} ^{\frac{3}{\gamma}}
 $$
 for $3< \gamma \leq \infty$.
Combining  (\ref{estfori1}) and (\ref{estfori2}) with
(\ref{eqfordu}), we obtain
 \bb\label{estfordu2}
  \frac{d}{dt} \|Du\|_{L^2} ^2 + \|D^2 u \|_{L^2} ^2 \leq C
  \|v\times \xi \|_{L^\gamma} ^{\frac{2\gamma}{\gamma -3}} \|Du
  \|_{L^2} ^2 +C \|v\|_{L^2}^2+C\|h\|_{L^2} ^2.
  \ee
By Gronwall's lemma we have
 \bq\label{gron}
  \lefteqn{\|Du (\cdot, t_0)\|_{L^2}
^2 +\nu \int_{t_0 -r_2^2} ^{t_0} \|D^2 u(\cdot,
 t)\|_{L^2}^2 dt}\hspace{.0in}\n \\
&& \leq \|Du (\cdot, t_0-r_2 ^2)\|_{L^2
 }^2\exp\left(C
 \int_{t_0-r_2^2} ^{t_0}
  \|v\times \xi (\cdot,t)\|_{L^\gamma } ^{\frac{2\gamma}{\gamma -3}
  }dt\right)\n \\
&&\qquad+C \int_{t_0-r_2^2} ^{t_0}\|h(\cdot,t)\|_{L^2 }^2 dt +
C\int_{t_0-r_2^2} ^{t_0}\|v(\cdot,t)\|_{L^2 }^2 dt.
  \eq
  Since $v\times \xi \in
L^{\gamma, \alpha}_{x,t} (Q_{z_0, r_2} )$ with $3/\gamma +2/\alpha
\leq 1$ and $\gamma > 3$, we estimate
 \bb\label{gron1}
\int_{t_0-r_2^2} ^{t_0}
  \|v\times \xi (\cdot,t)\|_{L^\gamma } ^{\frac{2\gamma}{\gamma -3}
  }dt\leq \|v\times \xi \|_{L^{\gamma, \alpha}_{x,t} (B_2 \times
  (t_0-r_2^2 ,t_0 ))} ^{\frac{\gamma}{\gamma-3}} r_2
  ^{\frac{2\gamma}{\gamma-3}
  (1-\frac{3}{\gamma}-\frac{2}{\alpha})} <\infty.
 \ee
From (\ref{gron})  and (\ref{gron1}) we find that $Du\in L^{2,
\infty}_{x,t} (Q_{z_0 ,r_2}) $, and hence $Dv\in L^{2, \infty}_{x,t}
(Q_{z_0 ,r_1})$. Similarly to the previous case, we
conclude that $z_0$ is a regular point for $v$.\\
Now,  we assume (iii) of the theorem holds true, and set $\eta
=v/|v|$. Then, we estimate in the preliminary step
 \bq\label{preest}
 \lefteqn{ |(\varphi v\times \o , D^2 u)_{L^2} |\leq \int_{B_2 } |\o \times \eta
  | |\varphi v| |D^2 u |dx}\hspace{.1in}\n \\
  &&=\int_{B_2 } |\o \times \eta | |u +V| |D^2 u |dx\n \\
  &&\leq \int_{B_2 } |\o \times \eta
  | |u | |D^2 u |dx
  +\int_{B_2 } |\o\times \eta || V| |D^2 u |dx\n \\
 &&\leq  \|\o \times \eta \|_{L^\gamma}
  \|u\|_{L^{\frac{2\gamma}{\gamma -2}}} \|D^2 u\|_{L^2}
  +\|V\|_{L^\infty}\|\o \|_{L^2} \|D^2 u\|_{L^2}\n \\
  &&=J_1+J_2.
  \eq
  The estimate of $J_2$ is simple as follows.
  \bb\label{estforj2}
  J_2 \leq C \|Dv\|_{L^2} ^2 + \frac{\nu}{8} \|D^2 u\|_{L^2}^2.
  \ee
For $2<\gamma \leq  3$, we have $6\leq
 \frac{2\gamma}{\gamma -2} <\infty $, and the following interpolation
 inequality is valid
 \bb\label{interpol1}
 \|u\|_{L^\frac{2\gamma}{\gamma -2}} \leq C \|u\|_{L^6}
 ^{2-\frac{3}{\gamma}} \|Du\|_{L^6} ^{-1+\frac{3}{\gamma}}\leq
 C \|Du\|_{L^2}
 ^{2-\frac{3}{\gamma}} \|D^2 v\|_{L^2} ^{-1+\frac{3}{\gamma}}.
 \ee
 If $\gamma=2$, then, instead,  we  use the inequality
 \bb\label{interpol1a}
\|u\|_{L^\infty} \leq C \|u\|_{L^6}^{\frac12}\|Du\|_{L^6} ^{\frac12}
\leq C \|Du\|_{L^2} ^{\frac12} \|D^2 u \|_{L^2} ^{\frac12}.
 \ee
 Substituting (\ref{interpol1}) or (\ref{interpol1a}) into
(\ref{preest}), we obtain
 \bq\label{estforj1}
J_1 &\leq& C\|\o \times \eta \|_{L^\gamma}\|Du\|_{L^2}
 ^{2-\frac{3}{\gamma}}\|D^2 u\|_{L^2} ^{\frac{3}{\gamma}}\n \\
&\leq& C\|\o \times \eta \|_{L^\gamma} ^{\frac{2\gamma}{2\gamma -3}}
\|Du\|_{L^2} ^2 +\frac{1}{8} \|D^2 u\|_{L^2} ^2
 \eq
 for $2\leq \gamma \leq 3$.
Combining (\ref{estforj1}) and (\ref{estforj2}), we have
 $$
    |(\varphi v\times\o , D^2 u)_{L^2}| \leq C\|\o \times \eta \|_{L^\gamma} ^{\frac{2\gamma}{2\gamma -3}}
\|Du\|_{L^2} ^2 +\frac{1}{4} \|D^2 u\|_{L^2} ^2 +C\|Dv\|_{L^2}^2.
$$
 Hence, from (\ref{preest}) we derive
 $$
  \frac{d}{dt} \|Du\|_{L^2} ^2 + \|D^2 u \|_{L^2} ^2 \leq C \|\o
  \times \eta \|_{L^\gamma} ^{\frac{2\gamma}{2\gamma -3}} \|Du
  \|_{L^2} ^2 +C \|h\|_{L^2} ^2+ \|Dv\|_{L^2}^2 .
  $$
  By Gronwall's lemma we have
   \bq\label{last1}
   \lefteqn{\|Du (\cdot, t_0)\|_{L^2}
^2 + \int_{t_0 -r_2^2} ^{t_0} \|D^2 u(\cdot,
 t)\|_{L^2}^2 dt}\hspace{.0in}\n \\
&& \leq \|Du (\cdot, t_0-r_2^2)\|_{L^2
 }^2\exp\left(C
 \int_{t_0-r_2^2} ^{t_0}
 \|\o
  \times \eta \|_{L^\gamma} ^{\frac{2\gamma}{2\gamma -3}}dt\right)\n \\
&&\qquad+C \int_{t_0-r_2^2} ^{t_0}\|Dv(\cdot,t)\|_{L^2 }^2 dt+C
\int_{t_0-r_2^2} ^{t_0}\|h(\cdot,t)\|_{L^2 }^2 dt.
 \eq
We  observe that
 \bb\label{last2}
 C\int_{t_0-r_2^2} ^{t_0}\|Dv(\cdot,t)\|_{L^2 }^2 dt+C
\int_{t_0-r_2^2} ^{t_0}\|h(\cdot,t)\|_{L^2 }^2 dt <\infty.
 \ee
Since $\o \times \eta \in L_{x,t}^{\gamma, \alpha}
  (Q_{z_0,r_2})$ with $3/\gamma +2/\alpha \leq 2$ and $2\leq \gamma \leq 3$ by hypothesis,
  we have
  \bb\label{finalest1}
\int_{t_0-r_2^2} ^{t_0}
  \|\o\times \eta(\cdot,t)\|_{L^\gamma } ^{\frac{2\gamma}{2\gamma -3}
  }dt\leq \|\o\times \eta \|_{L^{\gamma, \alpha}_{x,t} (B_2 \times
  (t_0-r_2^2 ,t_0 ))} ^{\frac{2\gamma}{2\gamma-3}} r_2
  ^{\frac{4\gamma}{2\gamma-3}
  (2-\frac{3}{\gamma}-\frac{2}{\alpha})} <\infty.
  \ee
From (\ref{last1})-(\ref{finalest1}) we find that $Du\in L^{2,
\infty}_{x,t} (Q_{z_0 ,r_2}) $,  and therefore $Dv\in L^{2,
\infty}_{x,t}
(Q_{z_0 ,r_1})$, and similar conclusion to the above case is obtained.\\
 Next, we assume that $z_0$ is a
singular point for which $t_0$ is not an epoch of possible
irregularity. Then, there exists a  time $t^*\in (t_0 -r^2, t_0)$
and $0<\tilde{r}_1<\tilde{r}_2 <r$ such that $v$ is regular on
$(B_{x_0,\tilde{ r}_2}\setminus B_{x_0,\tilde{ r}_1})\times[t^*,
t_0]$. This is due to that fact that the one dimensional Hausdorff
measure of the set of all possible singular space-time points is
equal to zero. We claim $v$ is regular on $B_{x_0,\tilde{ r}_1}
\times [t^*, t_0]$. Suppose not, then there exists another time
$s\in [t^*, t_0]$ such that the weak solution is regular on
$B_{x_0,\tilde{ r}_1} \times[t^*, s)$, and singularity occurs at
$(y,s)\in B_{x_0,\tilde{ r}_1} \times \{s\}$. We can repeat the
above argument for the parabolic neighborhoods of $(y,s)$ to
conclude that $(y,s)$ is actually a regular point. Hence, there
exists no space-time point of singularity in $ B_{x_0,\tilde{ r}_1}
\times [t^*, t_0]$, and we are reduced to the already considered
case that $t_0$ is an epoch of possible irregularity. This completes
the proof. $\square$


\begin{thebibliography}{1}
\bibitem{bei1} H. Beir\~{a}o da Veiga, {\it Vorticity and
smoothness in incompressible viscous flows,} in ``Wave Phenomena and
Asymptotic Analysis", RIMS, Kokyuroku 1315 (2003), pp. 37--45.
\bibitem{bei2} H. Beir\~{a}o da Veiga, {\it Concerning the
regularity problem for the solutions of the Navier-Stokes
equations,} C. R. Acad. Sci. Paris, Ser. I. Math. {\bf 321}, (1995),
pp. 405--408.
\bibitem{bei3} H. Beir\~{a}o da Veiga and L. C. Berselli,
{\it On the regularizing effect of the vorticity direction in
incompressible viscous flows,} Diff. Int. Eqns, {\bf 15}, (2002),
pp. 345--356.
\bibitem{caf} L. Caffarelli, R. Kohn, and L. Nirenberg, {\it Partial
regularity of suitable weak solutions of the Navier-Stokes
equations}, Comm. Pure Appl. Math., {\bf 35}, (1982), pp. 771--831.
\bibitem{cha1} D. Chae, {\it On the regularity conditions for the
Navier-Stokes and related equations},  Revista Mat. Iberoamericana.,
{\bf 23}, no. 1, (2007), pp. 371--384.
\bibitem{cha2} D. Chae and H.-J. Choe, {\it Regularity of
Solutions to the Navier-Stokes equations}, Electronic J. Diff.
Eqns., (1999), 7pp.(electronic).
\bibitem{cha3} D. Chae, K. Kang and J. Lee,
   {\it On the interior regularity of suitable weak solutions to
   the Navier-Stokes equations}, to appear in Comm. P.D.E.
\bibitem{con1} P. Constantin and C. Fefferman, {\it
Direction of vorticity and the problem of global regularity for the
Navier-Stokes equations,} Indiana Univ. Math. J., {\bf 42}, (1993),
pp. 775--789.
\bibitem{con2}P. Constantin and A. Majda, {\it The Beltrami spectrum
for incompressible fluid flows,} Comm. Math. Phys., {\bf 115},
(1988), pp. 435--456.
\bibitem{esc} L. Escauriaza, G. Seregin and V. Sverak,
{\it $L^{3,\infty}$-solutions of Navier-Stokes equations and
backward uniqueness,} Russian Math. Surveys, {\bf 58}, (2003), pp.
211--250.
\bibitem{fab} E. Fabes, B. Jones and N. Riviere, {\it The
initial value problem for the Navier-Stokes equations with data in
$L^p$}, Arch. Rat. Mech. Anal., {\bf 45}, (1972), pp. 222--248.
\bibitem{foi} C. Foias and R. Temam, {\it Some analytic and
geometric properties of the solutions of the evolution Navier-Stokes
equations}, J. Math. Pures Appl., {\bf 58}, (1979), pp. 339--368.
\bibitem{gal} G. P. Galdi, {An introduction to the Mathematical
theory of the Navier-Stokes initial-boundary value problem}, G. P.
Galdi, J. Heywood, and R. Rademacher(editors), Fundamental
directions in Mathematical Fluid Mechanics, Advances in Mathematical
Fluid Mechanics, Vol. 1, Birkh\"{a}user-Verlag, Basel, 2000.
\bibitem{hop} E. Hopf, {\"{U}ber die Anfangswertaufgabe f\"{u}r die
hydrodynamischen Grundgleichungen}, Math. Nachr., {\bf 4}, (1951),
pp.213--231.
\bibitem{lad} O. A. Ladyzhenskaya, {\it On the uniqueness and
smoothness of generalized solutions of the Navier-Stokes equations},
Zapiski Scient. Sem. LOMI, {\bf 5}, (1967), pp. 169--185.
\bibitem{ler} J. Leray, {\it Essai sur le mouvement d'un fluide
visqueux emplissant l'espace,} Acta Math., {\bf 63}, (1934),
pp.193--248.
\bibitem{lin} F. Lin, {\it A new proof of the
Caffarelli-Kohn-Nirenberg Theorem,} Comm. Pure Appl. Math., {\bf
51}, (1998), pp. 241--257.
\bibitem{neu} J. Neustupa and P. Penel, {\it Regularity of a
suitable weak solution to the Navier-Stokes equations as a
consequence of a regularity of one velocity component}, in: H.
Beir\~{a}o da Veiga, A. Sequeira and J. Videman (editors), Nonlinear
Applied Analysis, Plenum Press, New york, 1999, pp. 391--402.
\bibitem{ohy} T. Ohyama, {\it Interior regularity of weak
solutions to the Navier-Stokes equation}, Proc. Japan Acad., {\bf
36}, (1960), pp. 273--277.
\bibitem{pel} R. Pelz, V. Yakhot, S. Orszag, L. Shtilman and E.
Levich, {\it Velocity-vorticity patterns in turbulent flows,} Phys.
Rev. Lett., {\bf 54}, (1985), pp. 2505--2509.
\bibitem{pro} G. Prodi, {\it Un teorama di unicita per le equazioni di Navier-Stokes,}
  Annali di Mat., {\bf 48}, (1959), pp. 173--182.
\bibitem{sch} V. Scheffer, {\it Hausdorff measure and the
Navier-Stokes equations}, Comm. Math. Phys., {\bf 55}, (1977), pp.
97--112.
\bibitem{ser} J. Serrin, {\it On the interior regularity of
weak solutions of the Navier-Stokes equations},
 Arch. Rational Mech. Anal.,
 \bibitem{str} M. Struwe, {\it On partial regularity results for
the Navier-Stokes equations}, Comm. Pure Appl. Math., {\bf 41},
(1988), pp. 437--458.
\bibitem{tak} S. Takahashi, {\it On interior regularity
criteria for weak solutions of the Navier-Stokes equations,}
Manuscripta Math., {\bf 69}, (1990), pp. 237--254.
\bibitem{tia}G. Tian and Z. Xin, {\it Gradient estimation on Navier-Stokes
equations},  Comm. Anal. Geom., {\bf 7},  no. 2, (1999), pp.
221--257.
\bibitem{wol}J. Wolf, {\it The Caffarelli-Kohn-Nirenberg theorem-a
direct proof by Campanato's method,} preprint, (2007).
  \end{thebibliography}
\end{document}